\documentclass{article}
\usepackage{latexsym}
\usepackage{graphics}
\usepackage{amsfonts}
\newtheorem{thm}{Theorem}

\newtheorem{prop}{Proposition}

\newtheorem{defn}{Definition}
\newtheorem{cor}{Corollary}
\newtheorem{quest}{Question}

\begin{document}

\title{Intersection Graphs for String Links}
\author{Blake Mellor\\
			Mathematics Department\\
			Loyola Marymount University\\
			Los Angeles, CA  90045-2659\\
   {\it  bmellor@lmu.edu}}
\date{}
\maketitle

\begin{abstract}
We extend the notion of intersection graphs for chord diagrams in the theory of finite type
knot invariants to chord diagrams for string links.  We use our definition
to develop weight systems for string links via the adjacency matrix of the intersection graphs,
and show that these weight systems are related to the weight systems induced by the
Conway and Homfly polynomials.
\end{abstract}
\tableofcontents

\section{Introduction} \label{S:intro}

The theory of finite type invariants allows us to interpret many important knot and link
invariants in purely combinatorial terms, as functionals on spaces of chord diagrams.  For
knots, there is an obvious intersection graph associated with these diagrams, first studied
(in this context) by Chmutov, Duzhin and Lando \cite{cdl}.  In many cases, these graphs
contain all the relevant information for the functionals coming from knot invariants
\cite{cdl, me}.  Some of the most important knot invariants, such as the Conway, Jones,
Homfly and Kauffman polynomials, can be interpreted in terms of these intersection graphs
\cite{bg, me2}.

However, it is not obvious how to extend the idea of the intersection graph to finite type
invariants of other objects, such as braids, links and string links.  The goal of this paper is
to introduce a reasonable definition of the intersection graph for chord diagrams associated
with string links.  As evidence that this is the "right" definition, we use these intersection
graphs to construct weight systems for string links, and show that these weight systems are
related to (though weaker than) the weight systems arising from the Conway and Homfly
polynomials.  In future work we will show that in some cases the intersection graph contains all
the relevant information in the string link chord diagram \cite{me3, me4}, and that it can be
used to give a new interpretation of Milnor's homotopy invariants \cite{me5}.

In section~\ref{S:prelim} we will review the necessary background,
including finite type invariants, chord diagrams, intersection graphs
and Lando's graph bialgebra \cite{la}.  In section~\ref{S:IGstring} we
will define intersection graphs for string links, and provide the
space of graphs with a bialgebra structure similar to Lando's.  In section~\ref{S:2term} we
look at chord diagrams modulo the 2-term relations introduced by Bar-Natan and Garoufalides
\cite{bg}, and use the adjacency matrices of their intersection graphs to construct weight
systems related to those arising from the Conway and Homfly polynomials.  Finally, in
section~\ref{S:questions} we pose some questions for further research.\\
\\
\noindent{\it Acknowledgement:}  The author thanks Loyola Marymount University
for supporting this work via a Summer Research Grant in 2003.

\section{Preliminaries} \label{S:prelim}

\subsection{Knots, Links and String Links} \label{SS:ksl}

An {\it (oriented) knot} is an embedding of the (oriented) circle $S^1$ into
the 3-sphere $S^3$.  A {\it knot invariant} is a map from these embeddings to
some set which is invariant under isotopy of the embedding.  We
will also consider invariants of {\it regular} isotopy, where the
isotopy preserves the {\it framing} of the knot (i.e. a chosen
section of the normal bundle of the knot in $S^3$).  A {\it link} with $k$
components is simply an embedding of the disjoint union of $k$ copies of $S^1$
into $S^3$; each component is a knot.  A {\it string link} can be defined as
follows:

\begin{defn} (Habegger and Lin \cite{hl})
Let D be the unit disk in the plane and let I = [0,1] be the unit interval. 
Choose k points $p_1,..., p_k$ in the interior of D, aligned in order along
the the x-axis.  A {\bf string link} $\sigma$ of k components is a smooth
proper imbedding of k disjoint copies of I into $D \times I$:
$$\sigma:\ \bigsqcup_{i=1}^k{I_i} \rightarrow D \times I$$
such that $\sigma|_{I_i}(0) = p_i \times 0$ and $\sigma|_{I_i}(1) = p_i
\times 1$. The image of $I_i$ is called the ith string of the string link
$\sigma$.
\end{defn}

Link and string link invariants are defined in the same way as for knots.  Note that any string
link can be closed up to a link in a unique way by joining the top and bottom of each component
by an arc which lies outside of $D \times I$ (and which is unlinked with the other such arcs). 
So any link invariant gives rise to a string link invariant by evaluating it on the closure of
the string link.

\subsection{Finite Type Invariants} \label{SS:finitetype}

Our treatment of finite type invariants will follow the combinatorial approach of Birman
and Lin~\cite{bl}.  We will give a brief
overview of this combinatorial theory, and its natural extensions to links
and string links.  For more details, see Bar-Natan~\cite{bn}.

We first note that we can extend any link invariant to an invariant of {\it
singular} links, where a singular link is an immersion of
a disjoint union of copies of $S^1$ in 3-space which is an embedding except for a finite
number of isolated double points. Given a link invariant $v$, we extend it
via the relation:

$$\includegraphics{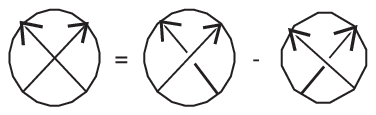}$$

An invariant $v$ of singular links is then said to be of {\it
finite type}, specifically of {\it type n}, if $v$ is zero on any
link with more than $n$ double points (where $n$ is a finite
nonnegative integer). We denote by $V_n$ the vector space over ${\mathbb
C}$ generated by (framing-independent) finite type invariants of
type $n$.  We can completely understand the space of finite type
invariants by understanding all of the vector spaces $V_n/V_{n-1}$.  An
element of this vector space is completely determined by its behavior on
links with exactly $n$ singular points.  In addition, since such an
element is zero on links with more than $n$ singular points, any other
(non-singular) crossing of the knot can be changed without
affecting the value of the invariant.  This means that elements of
$V_n/V_{n-1}$ can be viewed as functionals on the space of {\it
chord diagrams}:

\begin{defn}
A {\bf chord diagram of degree n with k components} is a disjoint union of k oriented
circles, together with $n$ chords (line segments with endpoints on the circles), such that
all of the $2n$ endpoints of the chords are distinct.  The circles represent a link of k
components and the endpoints of a chord represent 2 points identified by the immersion of
this link into 3-space.  The diagram is determined by the orders of the endpoints on
each component.  For string links, chord diagrams involve oriented line segments
rather than circles.
\end{defn}

Functionals on the space of chord diagrams which are derived from
finite type link invariants will satisfy certain relations:
\begin{defn}
A {\bf weight system of degree n} is a linear functional $W$ on
the space of chord diagrams of degree $n$ (with values in an
associative commutative ring ${\bf K}$ with unity) which satisfies
the 1-term and 4-term relations, shown in Figure~\ref{F:4-term}.  For knots,
the three arcs all belong to the same circle; for links (resp. string links),
they may belong to the same or different circles (resp. line segments).
    \begin{figure}
    (1-term relation) $$\includegraphics{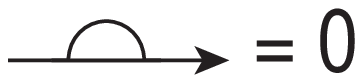}$$
    (4-term relation) $$\includegraphics{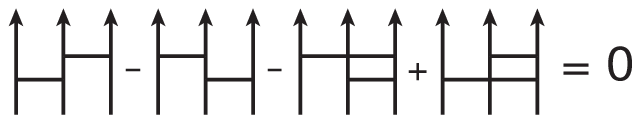}$$
    \caption{The 1-term and 4-term relations.  No other chords have endpoints
    on the arcs shown.  In the 4-term relations, all other chords of the four
    diagrams are the same.} \label{F:4-term}
    \end{figure}
\end{defn}

The natural map from elements of $V_n/V_{n-1}$ to functionals on chord
diagrams is a homomorphism into the space of weight systems
\cite{bn,bl,st,va}.  Kontsevich proved the much more difficult fact that
these spaces are isomorphic \cite{bn,ko} (the inverse map is the famous
{\it Kontsevich integral}).  For convenience, we take the dual approach,
and simply study the space of chord diagrams of degree $n$ modulo the
1-term and 4-term relations. The 1-term relation is occasionally referred
to as the "framing-independence" relation, because it arises from the
framing-independence of the invariants in $V_n$ (essentially, from
the first Reidemeister move). Since most of the interesting
structure of the vector spaces arises from the 4-term relation, it
is common to look at the more general setting of invariants of
regular isotopy, and consider the vector space $A_n^k$ of chord
diagrams of degree $n$ on links with $k$ components modulo the 4-term relation alone. 
Linear functionals on $A_n^k$ are called {\it regular weight systems} of degree $n$. 
Similarly, we define the vector space $B_n^k$ for chord diagrams of degree $n$ on {\it
string links} of $k$ components.

It is useful to combine all of these spaces into graded modules
$A^k = \bigoplus_{n\geq 1}A_n^k$ and $B^k = \bigoplus_{n\geq 1}B_n^k$ via direct sum. 
For $k = 1$, we can define a product on $A^1$ via connect sum \cite{bn}; however,
this does not extend to $k > 1$.  But we {\it can} give the module
$B^k$ a bialgebra (or Hopf algebra) structure for any $k$ by defining an appropriate
product and co-product:
\begin{itemize}
    \item  We define the (noncommutative) product $D_1 \cdot D_2$ of two chord diagrams
$D_1$ and $D_2$ as the result of placing $D_2$ on top of $D_1$ (joining the components
so the orientations agree), as shown below:
$$\includegraphics{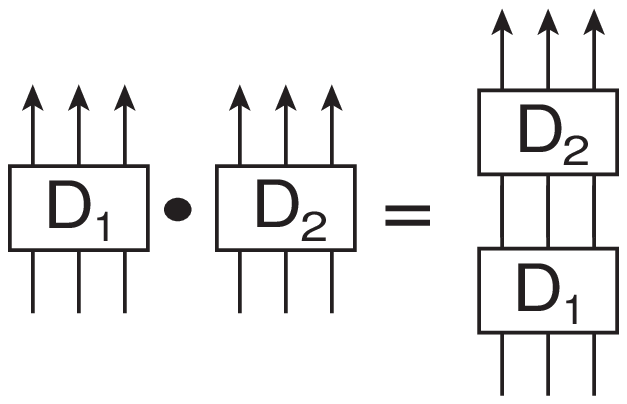}$$
    \item  We define the co-product $\Delta(D)$ of a chord diagram $D$ as
follows:
$$\Delta(D) = {\sum_J D_J' \otimes D_J''}$$
where $J$ is a subset of the set of chords of $D$, $D_J'$ is $D$
with all the chords in $J$ removed, and $D_J''$ is $D$ with all
the chords {\it not} in $J$ removed.
$$\includegraphics{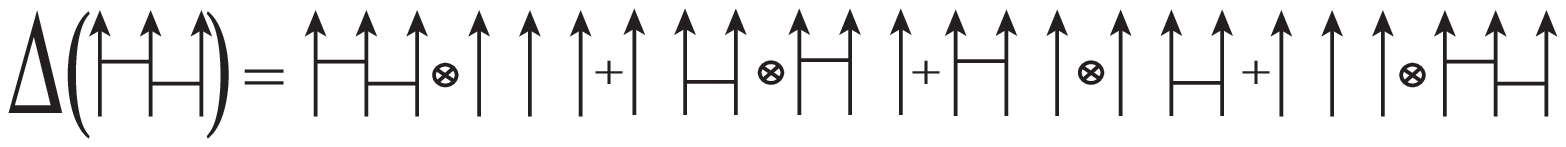}$$
\end{itemize}
It is easy to check the compatibility condition $\Delta(D_1\cdot
D_2) = \Delta (D_1)\cdot\Delta(D_2)$.

\section{Intersection Graphs for String Links} \label{S:IGstring}

For knots, there is a very natural notion of intersection among the chords of
a chord diagram - two chords intersect if their endpoints alternate around
the bounding circle.  The intersection graph then has a vertex for each
chord, and connects two vertices by an edge if the corresponding chords
intersect.  The difficulty with extending this idea to links and string links
is how to deal with chords which have their endpoints on different components
of the diagram - specifically, if they each have one endpoint on a given
component, and their other endpoints are on different components.  In what
sense, if any, do these chords intersect?

\subsection{Definition} \label{SS:def}

The essential value of the intersection graph for knots is that it can detect
when the order of two endpoints for different chords along the bounding circle
is switched (as happens in the 4-term relation), since this changes the pair
of chords from intersecting to non-intersecting or vice-versa.  To usefully
extend the notion of intersection graphs to links or string links, we need to retain this
ability.  For string links, the existence of a "bottom" and "top" for each component allows
us to give a linear (rather than cyclic) ordering to the endpoints of the chords on each
component, and so the notion of one endpoint being "below" another is well-defined.

\begin{defn}
Let $D$ be a string link chord diagram with $k$ components (oriented line
segments, colored from 1 to $k$) and $n$ chords.  The {\it intersection graph}
$\Gamma(D)$ is the labeled, directed multigraph such that:
\begin{itemize}
    \item $\Gamma(D)$ has a vertex for each chord of $D$.  Each vertex is labeled
by an unordered pair $\{i,j\}$, where $i$ and $j$ are the labels of the components
on which the endpoints of the chord lie.
    \item There is a directed edge from a vertex $v_1$ to a vertex $v_2$ for each
pair $(e_1, e_2)$ where $e_1$ is an endpoint of the chord associated to $v_1$,
$e_2$ is an endpoint of the chord associated to $v_2$, $e_1$ and $e_2$ lie on the
same component of $D$, and the orientation of the component runs from $e_1$ to
$e_2$ (so if the components are all oriented upwards, $e_1$ is below $e_2$).  We
count these edges "mod 2", meaning that if two vertices are connected by two
directed edges with the same direction, they cancel each other.  If two vertices
are connected by a directed edge in each direction, we will simply connect them
by an undirected edge.
\end{itemize}
\end{defn}

\begin{figure}
    $$\includegraphics{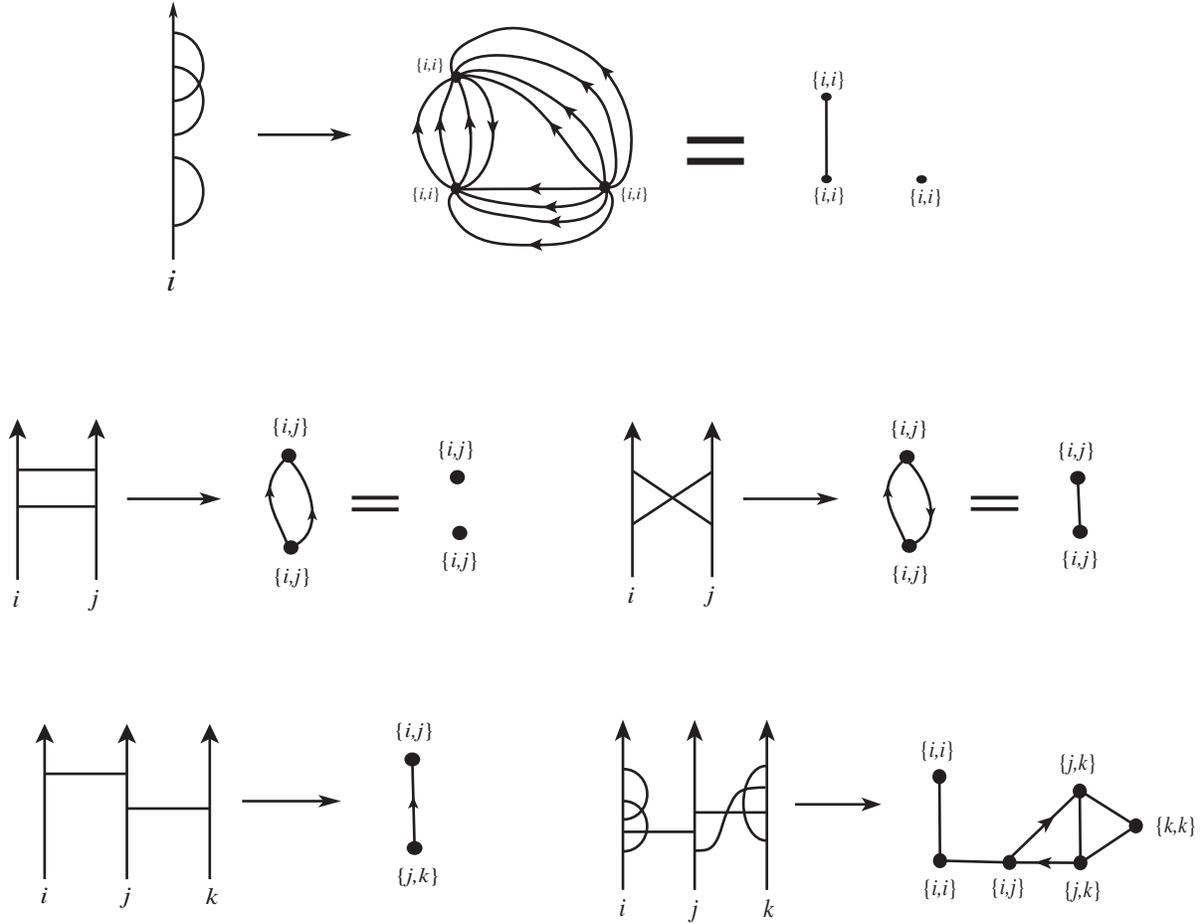}$$
    \caption{Examples of intersection graphs for string links} \label{F:IGstring}
    \end{figure}

Examples of chord diagrams and their associated intersection graphs are given in
Figure~\ref{F:IGstring}.
Note that when the two chords have both endpoints on the same component $i$, our
definition of intersection graph corresponds to the usual intersection graph for
knots.  Our definition also matches our intuition in the case of chord diagrams
of two components, as shown in Figure~\ref{F:IGstring}.

Note also that the total number of directed edges between a vertex $v$ labeled
$\{i,j\}$ and a vertex $w$ labeled $\{l,m\}$ is given by the sum of the number
of occurrences of $i$ in $\{l,m\}$ and the number of occurrences of $j$ in
$\{l,m\}$.  In particular, if a vertex $v$ has a label $\{i,i\}$, this number
will be even (0, 2 or 4).  Since we count directed edges modulo 2, this implies
there is an (uncancelled) directed edge from $v$ to another vertex $w$ if and only if
there is also an (uncancelled) directed edge from $w$ to $v$.  We will say that labeled
directed multigraphs which have this property are {\it semisymmetric}.

\begin{defn} \label{D:semisym}
A directed multigraph G, with each vertex labeled by a pair \{i,j\}, is {\bf
semisymmetric} if for every vertex v labeled \{i,i\}, and any other vertex w, there is a
directed edge from v to w if and only if there is a directed edge from w to v.
\end{defn}

\subsection{Graph bialgebra for String Links} \label{SS:bialgebra}

Following Lando's work for knots \cite{la}, we can define a bialgebra structure
on the space of intersection graphs for string links so that $\Gamma$ becomes a
bialgebra homomorphism from the space of chord diagrams for string links to the space of
intersection graphs.  The key is to define the analogue of the 4-term relation
for intersection graphs.

\begin{defn}
Consider the graded vector space (over ${\mathbb C}$) of formal
linear combinations of labeled semisymmetric directed multigraphs, with vertices
labeled by unordered pairs \{i,j\}, $1 \leq i,j \leq k$, graded by the number of
vertices in the graphs.  For any graph G and vertices A and B in V(G), with
labels \{i,j\} and
\{i,l\} respectively (j and l may be equal, or equal to i), we impose on the
vector space the relation:
$$G-G'_{AB}-\widetilde{G}_{AB}+\widetilde{G}'_{AB} = 0$$
Here $G'_{AB}$ is the result of complementing the edge AB in G
(i.e. adding a directed edge between the vertices in each direction, and then
cancelling "mod 2").  $\widetilde{G}_{AB}$ is the result of changing the
label on A to \{j,l\}, complementing the edge AB if $i \neq l$ (and leaving it
unchanged if i = l), and adding a directed edge from A to C (respectivley C to A)
for every vertex C in V(G) for which G has a directed edge from B to C
(respectively C to B), and then cancelling mod 2.  Finally,
$\widetilde{G}'_{AB}$ is the result of complementing the edge AB in
$\widetilde{G}_{AB}$.  Here is an example of such a relation:
$$\includegraphics{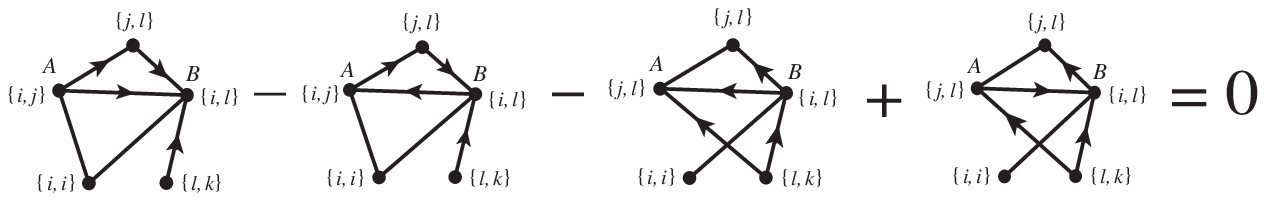}$$
The bialgebra $F$ is defined as this graded vector space, together
with a product and a coproduct.  The product is a map $\cdot: F \times F
\rightarrow F$, defined as follows.  Given graphs $G_1$ and $G_2$, $G_1\cdot G_2$
is the disjoint union of the graphs, together with a directed edge from $v_1 \in
V(G_1)$ to $v_2 \in V(G_2)$ for each color in the label \{i,j\} for $v_1$ which
is also in the label for $v_2$.  An example is shown below:
$$\includegraphics{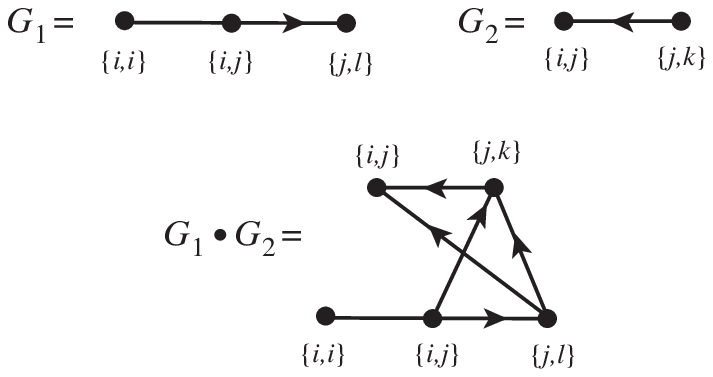}$$
The coproduct is a map $\mu: F \rightarrow F \otimes F$, defined as follows.  For
any graph G, and subset $J \subseteq V(G)$ of its vertices, let $G_J$ denote the
subgraph induced by $J$.  Then:
$$\mu(G) = \sum_{J \subseteq V(G)}{G_J \otimes G_{V(G)\backslash J}}$$
An example is shown below:
$$\includegraphics{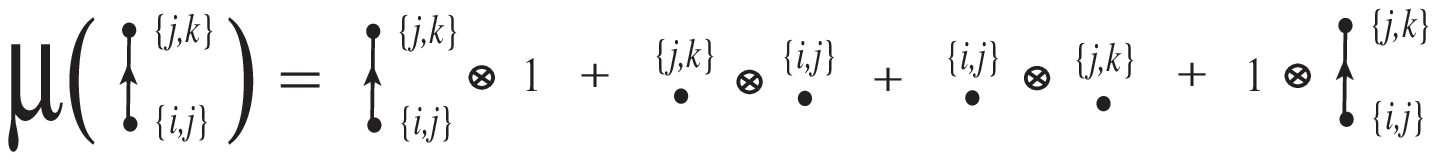}$$
\end{defn}

\begin{thm} \label{T:bialgebra}
The product and coproduct defined above induce the structure of a co-commutative
(but not commutative) bialgebra on the space F of labeled, directed multigraphs
modulo the 4-term relation.
\end{thm}
{\sc Proof:} It is easy to check that $F$ is a bialgebra (i.e. that $\mu(G_1\cdot
G_2) = \mu(G_1)\cdot\mu(G_2)$), and that it is co-commutative, but not
commutative.  We also need to check that the product and coproduct respect the
4-term relation, i.e. that the product of a graph and a 4-term relation yields a
4-term relation, and that the coproduct of a 4-term relation yields a sum of
4-term relations.  To begin with, consider a product
$(G-G'_{AB}-\widetilde{G}_{AB}+\widetilde{G}'_{AB})\cdot H =
G\cdot H-G'_{AB}\cdot H-\widetilde{G}_{AB}\cdot H+\widetilde{G}'_{AB}\cdot H$.  It
is clear that $G'_{AB}\cdot H = (G\cdot H)'_{AB}$, since complementing the edge
$AB$ is independent of any edges between $G$ and $H$.  Moreover,
$\widetilde{G}_{AB}\cdot H = \widetilde{G\cdot H}_{AB}$, since in each case there
is an directed edge from a vertex $v \in V(H)$ to $A$ for each occurrence of $j$
or $l$ in the label of $v$ (with cancellation mod 2).  It is then clear that
$\widetilde{G}'_{AB}\cdot H = \widetilde{G\cdot H}'_{AB}$, and so the product of
the 4-term relation and another graph is another 4-term relation.

For the coproduct, consider
$\mu(G-G'_{AB}-\widetilde{G}_{AB}+\widetilde{G}'_{AB})$.  The sum in the
coproduct splits into two groups - terms where both $A$ and $B$ belong to either
$J$ or its complement, and terms where one of the two vertices is in $J$ and the
other is in its complement.  The first group gives a sum of 4-term relations,
while the terms of the second group already sum to zero in $\mu(G-G'_{AB})$ and
$\mu(\widetilde{G}_{AB}-\widetilde{G}'_{AB})$.  This finishes the proof of the
theorem. $\Box$

It is now easy to check that $\Gamma$ is a bialgebra homomorphism from the
bialgebra $B^k$ of chord diagrams for string links with $k$ components modulo the 4-term
relation to the bialgebra $F^k$ of directed, labeled graphs with labels 1,..., $k$, since
the 4-term relation for graphs was defined to mimic the 4-term relation for the chord
diagrams.  As a result, any functional on $F^k$ will induce, by composition with $\Gamma$, a
regular weight system.  We call such a functional a {\it regular graph weight system}.

\section{2-term relations for String Links} \label{S:2term}

\begin{figure} [h]
    $$\includegraphics{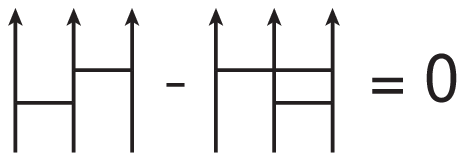}$$
    \caption{The 2-term relations for string links} \label{F:2-term}
    \end{figure}

Any particular weight system will satisfy relations in addition to the 1-term and 4-term
relations, and it can be useful to look at weight systems which lie in the subspaces
determined by these additional relations.  In particular, Bar-Natan and Garoufalides
\cite{bg} noted that the weight system associated with the Conway polynomial for knots
satisfies a set of {\it 2-term} relations which induce the 4-term relations.  We can extend
this notion to define a set of 2-term relations for string links, shown in
Figure~\ref{F:2-term}.
Clearly, these relations imply that the weight system satisfies the 4-term relation as
well.  As a result, the product and coproduct of section \ref{SS:finitetype} are still
well-defined. So we can give the vector space of chord diagrams for string links with $k$
components modulo the 2-term relations the structure of a bialgebra.  We will denote this
bialgebra (and the underlying vector space) by $D^k$.  There is a natural projection from
$B^k$ to $D^k$.

Similarly, we can look at the quotient space of $F$ by the 2-term relations
$$G-\widetilde{G}_{AB} = 0$$
Since $\widetilde{G}'_{AB} = \widetilde{G'}_{AB}$, the 2-term relations imply the
4-term relations, and the quotient space is a bialgebra $E$ with the product and
coproduct induced from $F$.  $\Gamma$ is then a bialgebra homomorphism from $D^k$
to $E^k$.

\subsection{Chord diagrams for String Links modulo 2-term relations} \label{SS:chord2term}

\begin{figure} [h]
    $$\includegraphics{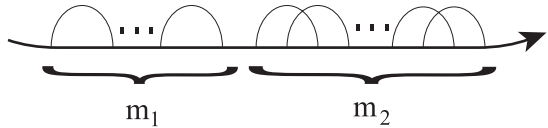}$$
    \caption{Example of an ($m_1, m_2$)-caravan} \label{F:caravan}
    \end{figure}

Bar-Natan and Garoufalides \cite{bg} analyzed the space of chord diagrams for knots
(equivalently, string links of one component) modulo the 2-term relations.  They found
that $D^1$ is generated (as a vector space) by {\it ($m_1,m_2$)-caravans} of $m_1$
"one-humped camels" (isolated chords which intersect no other chords) and $m_2$
"two-humped camels" (pairs of chords which intersect each other, but no other chords). An
example of such a caravan is shown in Figure~\ref{F:caravan}.

Our goal in this section is to find a similar "normal form" for string link chord
diagrams with more than one component; this will be useful in the remainder of the paper. 
We will prove the following theorem:

\begin{thm} \label{T:normalform}
Any connected string link chord diagram is equivalent, modulo the 2-term relations, to a
chord diagram with components numbered 1,..., n (possibly after renumbering) such that:
\begin{itemize}
				\item Every chord either has both endpoints on component 1, or endpoints on
components i and i+1 for some i.
				\item The chords with both endpoints on component 1 are arranged in a caravan. 
Moreover, all of their endpoints lie below the endpoints of any chords between components.
				\item There are at most 2 chords between components i and i+1, and if there are 2
they do not cross (so the endpoints of one chord lie above the endpoints of the other on
both components).
    \item On component i, the endpoints of chords connecting component i to component
i+1 lie below the endpoints of chords connecting component i to component i-1.
\end{itemize}
An example is shown in Figure~\ref{F:normal}.
    \begin{figure} [h]
    $$\includegraphics{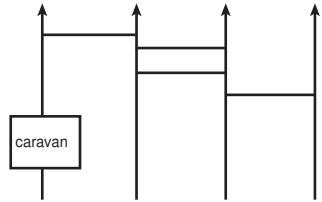}$$
    \caption{Example of a link chord diagram in normal form} \label{F:normal}
    \end{figure}
\end{thm}

If a diagram is not connected, we simply put each connected component into normal form.
 The proof proceeds in several steps.  We will describe an algorithm to put any connected
chord diagram into the normal form described above via 2-term relations.\\

\begin{figure} [h]
    $$\includegraphics{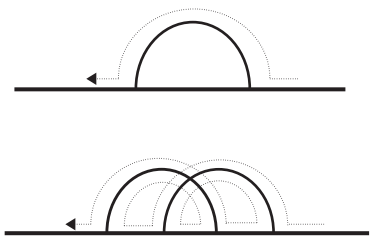}$$
    \caption{Factoring out one- and two-humped camels} \label{F:factor}
    \end{figure}

\begin{figure} [h]
    $$\includegraphics{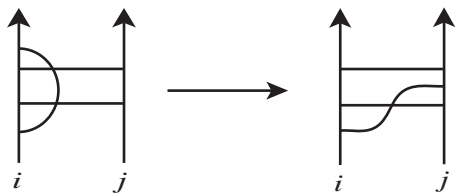}$$
    \caption{A move from Step 1} \label{F:step1}
    \end{figure}

\noindent{\sc Step 1:}  We first deal with chords which have both endpoints on the same
component.  To begin with, any pair of intersecting chords which have all four endpoints
on the same component can be slid to the bottom of the component, away from any other
chords, to form a two-humped camel as in \cite{bg}.  An isolated chord with both endpoints
on the same component can also be slid down to the bottom, forming a one-humped camel. 
See Figure~\ref{F:factor} for examples of these moves.
We are left with chords with both endpoints on the same component which only intersect
chords which lie between two different components.  By sliding this chord over one of the
chords it intersects, its endpoints are now on two different components, as shown in
Figure~\ref{F:step1}.

\begin{figure} [h]
    $$\includegraphics{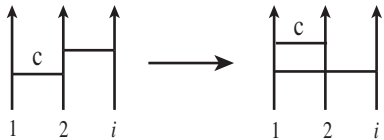}$$
    \caption{Moving an endpoint below $c$} \label{F:step2}
    \end{figure}

\noindent{\sc Step 2:}  Our next step is to put every chord connecting two different
components on a different horizontal "level"; in particular, we will remove intersections
between chords between the same two components.  Say that we have $n$ chords between
components.  Select a chord arbitrarily, which we will denote chord $c$.  Without loss of
generality, assume $c$ connects components 1 and 2.  We will move $c$ to the top level. 
As we move it up, if we come to a chord with an endpoint on components 1 or 2 which is
above $c$ we move that endpoint below $c$ via a 2-term relation, as shown in
Figure~\ref{F:step2}.
Ultimately, the endpoints of $c$ will lie above the endpoints of any other chord.  During
this process, some chords may have had both their endpoints moved to the same component -
remove or modify these chords as in Step 1 (notice that this only reduces the number of
chords between components from $n$).  We are left with {\it at most} $n-1$ chords
connecting components, all below $c$.  So we can continue the process inductively, ending
with all chords between components horizontal, at different levels (and with a caravan at
the bottom of each component).\\

\begin{figure} [h]
    $$\includegraphics{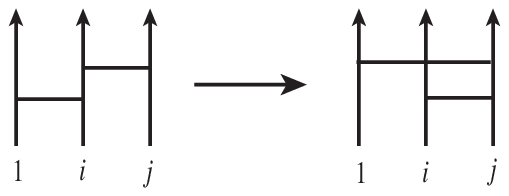}$$
    \caption{Step 3} \label{F:step3}
    \end{figure}

\noindent{\sc Step 3:}  Now we consider one component of the diagram, say component 1, and
move all the chords with an endpoint on component 1 to a higher level than all the other
chords in the diagram (leaving chords with both endpoints on the same component at the
bottom).  We will denote by an {\it (i,j)-chord} a chord with endpoints on components $i$
and $j$.  If a (1,$i$)-chord lies directly below an ($i,j$)-chord along component $i$, we
can use a 2-term relation to move the (1,$i$)-chord above the other chord, in the process
transforming it into a (1,$j$)-chord, as in Figure~\ref{F:step3}.
Since there are only a finite number of chords, we can successively move all chords with
one endpoint on component 1 above all other chords (while the other endpoint may move,
one will always remain on component 1).\\

\begin{figure} [h]
    $$\includegraphics{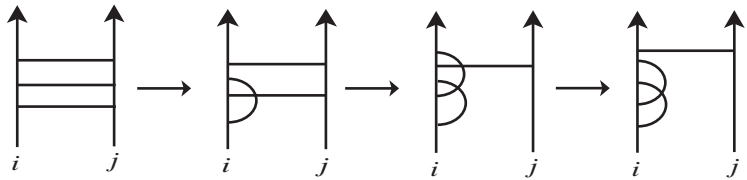}$$
    \caption{Factoring a two-humped camel from 3 parallel chords} \label{F:reduce}
    \end{figure}

\noindent{\sc Step 4:}  Let $c$ be the top chord at this point; without loss of generality,
$c$ is a (1,2)-chord.  Let $k_{1i}$ be the number of (1,$i$)-chords, and $k_1 =
\sum_i{k_{1i}}$.  Now move $c$ down component 1, stripping other chords off of component
1 by 2-term relations, transforming them from (1,$i$)-chords to (2,$i$)-chords, reversing
the move from Step 2 shown in Figure~\ref{F:step2}.  Continue this until $c$ encounters
another (1,2)-chord, $d$.  Now move $c$ back up as in Step 3 (see Figure \ref{F:step3});
when it is again at the top, it will be a (1,$i$)-chord for some $i$.  Repeat the process
for $d$, moving it up to just below $c$, and with all the other (1,2)-chords.  This
reduces $k_1$ to $k_{12}$, though the chords may no longer be (1,2)-chords.  We can
relabel the components so that $c$ is once again a (1,2)-chord and repeat the process. 
Each repetition reduces $k_1$, until we reach a point where $k_1 = k_{12}$ (i.e. all
chords with an endpoint on component 1 are (1,2)-chords), and the (1,2)-chords lie above
all other chords in the diagram.  Finally, we can reduce $k_{12}$ to 1 or 2 by noticing
that whenever $k_{12} \geq 3$, we can factor out a 2-humped camel and reduce $k_{12}$ by
2, as in Figure~\ref{F:reduce}.\\

\begin{figure} [h]
    $$\includegraphics{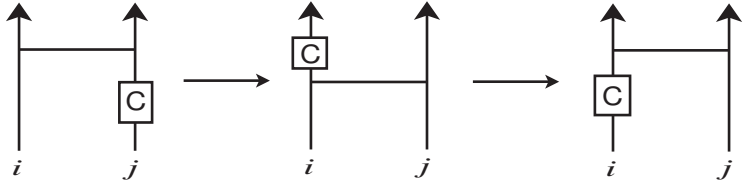}$$
    \caption{Moving a caravan} \label{F:step5}
    \end{figure}

\noindent{\sc Step 5:}  We can now repeat Steps 3 and 4 for component 2, and then for each
component in turn.  We are left with a diagram which is almost in normal form - the
final step is to move all of the caravans to component 1.  We can move a caravan from the
$i$th component to the $(i-1)$th component by sliding it over a chord between the two
components, as in Figure~\ref{F:step5}.
Any other chords between the two components can be slid over the caravan in turn by
2-term relations, leaving the caravan below the chords connecting components $i-1$ and
$i-2$.  Continuing in this way, the caravans can be moved to component 1, and finally
pushed to to bottom of component 1 as in step 1.  This completes the proof of Theorem
\ref{T:normalform}. $\Box$

\subsection{Adjacency Matrices for Intersection Graphs} \label{SS:adj}

In this section we will define the adjacency matrix for an intersection graph,
and show that its rank and determinant (over ${\bf Z}_2$) are regular graph weight
systems.

\begin{defn} \label{D:adj}
Given a directed graph G with n vertices, $\{v_1,...,v_n\}$, such that each
vertex $v_i$ has a label $\{a_i,b_i\}$, the {\bf adjacency matrix of G}, or
adj(G), is the $n\times n$ matrix defined by:
$$adj(G)_{ij}\ (i \neq j) = \left\{{\matrix{1\ if\ there\ is\ a\ directed\ edge\
from\ v_i\ to\ v_j\ in\ G \cr 0\ otherwise}}\right.$$
$$adj(G)_{ii} = \left\{{\matrix{1\ if\ a_i \neq b_i \cr 0\ if\ a_i = b_i}}\right.$$
\end{defn}

These matrices can be viewed as bilinear forms over ${\bf Z}_2$, but unlike the
case for knots \cite{me2}, they are not generally symmetric, so their
classification is much more difficult (and still, to my knowledge, an open
question).  We are interested in the notion of {\it congruence} of such
matrices.  We are particularly interested in the adjacency matrices of {\it
semisymmetric} graphs; which have the property (which we will also denote by
semisymmetry) that if the $i$th diagonal element is 0, then the $i$th row and
column are the same.

\begin{defn}
We say that two $n \times n$ matrices A and B are {\bf congruent} (denoted $A \cong B$) if
there is an invertible matrix P such that $A = PBP^{T}$.
\end{defn}

One of our main results is that the congruence class of the adjacency matrix
(over ${\bf Z}_2$) satisfies the 2-term relation, and so induces a regular graph
weight system.

\begin{thm}
For any graph G in E, $adj(G) \cong adj(\widetilde{G}_{AB})$.
\end{thm}
{\sc Proof:}  Reordering the vertices of $G$ changes $adj(G)$ by a congruence; in
this case, $P$ is the result of doing a corresponding reordering of the rows of
the identity matrix.  So we can assume that the first two rows and columns of
$adj(G)$ and $adj(\widetilde{G}_{AB})$ correspond to the vertices $A$ and $B$. 
The matrices are identical except for the first row and column.  Say that $A$
has label $\{c,a\}$ and $B$ has label $\{c,b\}$, so in $\widetilde{G}_{AB}$ $A$ has label
$\{a,b\}$.  An entry $adj(\widetilde{G}_{AB})_{1i}$ (where $i \neq 1,2$) is equal to
$adj(G)_{1i}+adj(G)_{2i}$ (mod 2), and similarly for the first column. 
$adj(\widetilde{G}_{AB})_{12} = adj(G)_{12}$ if and only if $b=c$, i.e. if
$adj(G)_{22}=0$, so $adj(\widetilde{G}_{AB})_{12} = adj(G)_{12}+adj(G)_{22}$. 
Finally, $adj(\widetilde{G}_{AB})_{11} = 0$ if $a=b$, and 1 otherwise.  We have
five cases.\\
\\
\noindent {\sc Case 1:} $a=b=c$.  In this case $adj(G)_{11} = adj(G)_{22} = 0$ and
$adj(G)_{12} + adj(G)_{21} = 4$ (due to the semisymmetry of the graph).  So, mod 2,
$adj(\widetilde{G}_{AB})_{11} = 0 = adj(G)_{11} + adj(G)_{22} + adj(G)_{12} +
adj(G)_{21}$.\\
\\
\noindent {\sc Case 2:} $a=b\neq c$.  In this case $adj(G)_{11} = adj(G)_{22} = 1$ and
$adj(G)_{12} + adj(G)_{21} = 4$.  So, mod 2, $adj(\widetilde{G}_{AB})_{11} = 0 =
adj(G)_{11} + adj(G)_{22} + adj(G)_{12} + adj(G)_{21}$.\\
\\
\noindent {\sc Case 3:} $a\neq b=c$.  In this case, $adj(G)_{11} = 1$ and $adj(G)_{22} = 0$,
and $adj(G)_{12} + adj(G)_{21} = 2$.  So, mod 2, $adj(\widetilde{G}_{AB})_{11} = 1 =
adj(G)_{11} + adj(G)_{22} + adj(G)_{12} + adj(G)_{21}$.\\
\\
\noindent {\sc Case 4:} $b\neq a=c$.  In this case, $adj(G)_{11} = 0$ and $adj(G)_{22} = 1$,
and $adj(G)_{12} + adj(G)_{21} = 2$.  So, mod 2, $adj(\widetilde{G}_{AB})_{11} = 1 =
adj(G)_{11} + adj(G)_{22} + adj(G)_{12} + adj(G)_{21}$.\\
\\
\noindent {\sc Case 5:} $a,b,c$ all different.  In this case $adj(G)_{11} = adj(G)_{22} = 1$
and $adj(G)_{12} + adj(G)_{21} = 1$.  So, mod 2, $adj(\widetilde{G}_{AB})_{11} = 1 =
adj(G)_{11} + adj(G)_{22} + adj(G)_{12} + adj(G)_{21}$.\\

We conclude that $adj(\widetilde{G}_{AB})_{11} = adj(G)_{11} + adj(G)_{22} + adj(G)_{12} +
adj(G)_{21}$.  Therefore $adj(G)$ and $adj(\widetilde{G}_{AB})$ are congruent by a matrix
$P$, where $P$ is the elementary matrix constructed from the identity matrix by adding the
second row to the first row.  In other words, $P$ is the same as the identity, except for
the 2 by 2 matrix in the upper left corner, which is $\left[{\matrix{1 & 1 \cr 0 &
1}}\right]$. $\Box$

\begin{cor}
The congruence class of the adjacency matrix of the intersection graph is
invariant under the 4-term relations.
\end{cor}

\begin{cor}
The determinant (mod 2) and rank (over ${\bf Z}_2$) of the adjacency matrix of the
intersection graph are regular graph weight systems.
\end{cor}
{\sc Proof:}  If $P$ is invertible over ${\bf Z}_2$, the rank of $PAP^T$ is the same
as the rank of $A$.  In addition, $det(P) = det(P^T) = 1$ (mod 2), so the determinant 
is also a congruence invariant mod 2.  $\Box$

\subsection{Relations to the Conway and Homfly polynomials} \label{SS:conhom}

For knots, the Conway and Homfly weight systems can be interpreted via intersection graphs
\cite{me2}.  To what extent can we do the same for string links?  In general, this is
still an open problem, but we can offer a few initial results.

\begin{figure} [h]
    $$\includegraphics{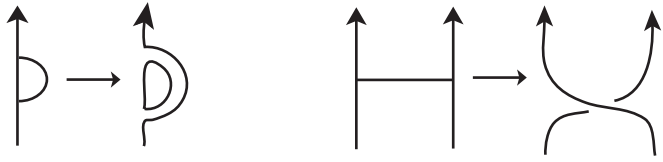}$$
    \caption{Surgery on a chord $v$} \label{F:surgery}
    \end{figure}

Let us recall the definitions of the Conway polynomial and weight system.  The Conway
polynomial $\Delta$ of a link is a power series $\Delta(L) = \sum_{n \ge 0}{a_n(L)z^n}$. 
It can be computed via the skein relation (where $L_+, L_-, L_0$ are as in
Figure~\ref{F:skein}):
    \begin{figure} [b]
    $$\includegraphics{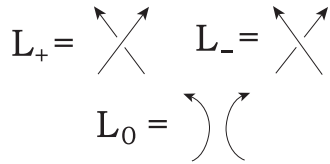}$$
    \caption{Diagrams of the skein relation.} \label{F:skein}
    \end{figure}
$$\Delta(L_+)-\Delta(L_-) = z\Delta(L_0)$$
$$\Delta(unlink\ of\ k\ components) = \left\{{\matrix{1\ if\ k=1 \cr 0\ if\ k>1}}\right.$$
We define the Conway polynomials of a string link as the Conway polynomial of the link
formed by closing the string link in the natural way.  The coefficient $a_n$ is a finite
type invariant of type $n$ \cite{bn,bl}, and therefore defines a weight system
$b_n$ of degree $n$.  The collection of all these weight systems is
called the Conway weight system, denoted $C$.  Consider a chord
diagram $D$, together with a chord $v$.  Let $D_v$ be the result
of {\it surgery on v}, i.e. replacing $v$ by an band which preserves the orientation of
the components, and then removing the interior of the band and the intervals where it is
attached to $D$, as shown in Figure~\ref{F:surgery}
(so $D_v$ may have multiple boundary circles).  The skein
relations for the Conway polynomial give rise to the following
relations for $C$:
$$C(D) = C(D_v)$$
$$C(unlink\ of\ k\ components) = \left\{{\matrix{1\ if\ k=1 \cr 0\ if\ k>1}}\right.$$

It is easy to show \cite{bg} that this weight system satisfies the 2-term relations of
section \ref{S:2term}.  Simply surger the two chords; the 2-term relation then says just
that one band can be "slid" over the other, which doesn't change the topology of the
diagram.  Bar-Natan and Garoufalides also showed that the weight system for the Conway
polynomial for knots is just the determinant of the adjacency matrix for the intersection
graph, mod 2.  The next theorem extends this result to string links with two components.

\begin{thm} \label{T:conway2}
For any chord diagram D on a string link with two components, $C(D) =
det(\Gamma(D))rank(\Gamma(D))$ (mod 2).
\end{thm}
{\sc Proof:}  Since both of these weight systems satisfy the
2-term relations, it suffices to show that they agree on diagrams in normal form (see
section \ref{SS:chord2term}). For string links with two components, normal form consists
of a caravan on the first component, and 0, 1, or 2 parallel chords between the two
components.  Say that the caravan has $m_1$ 1-humped camels and
$m_2$ 2-humped camels, and there are $l$ chords between the components ($l=0,1,2$). 
Then $adj(\Gamma(D)) \cong [1]^l \oplus [0]^{m_1} \oplus \left[{\matrix{0 & 1 \cr 1
& 0}}\right]^{m_2}$.  So $det(\Gamma(D)) = 1^{l+m_2}\cdot 0^{m_1}$
(mod 2) and $rank(\Gamma(D)) = l + 2m_2 \equiv l$ (mod 2).  The product of the rank and
determinant is 1 (mod 2) when $l = 1$ and $m_1 = 0$, and 0 otherwise.

On the other hand, if we surger all the chords of the diagram, and close the string link, we
obtain an unlink with $m_1 + 2$ components when $l=0,2$, and $m_1 + 1$ when $l=1$, which
means that $C(D) = 1$ only when $m_1=0$ and $l=1$.  So the two weight systems
agree.  $\Box$

\begin{figure} [h]
    $$\includegraphics{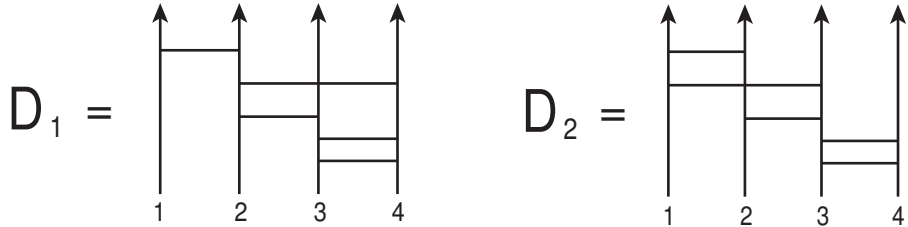}$$
    \caption{$C(D)$ is not determined by the adjacency matrix} \label{F:counterexample}
    \end{figure}

With more than two components, the situation is more complicated; in fact, the Conway
weight system is not always determined by the adjacency matrix of the intersection graph. 
Figure~\ref{F:counterexample} shows two chord diagrams $D_1$ and $D_2$ on 4 components
which illustrate this.  
It is easy to see that these intersection graphs for these diagrams differ only in
the labeling of the vertices, so:

$$adj(\Gamma(D_1)) = adj(\Gamma(D_2)) = \left[{\matrix{1 & 0 & 0 & 0 & 0 \cr 
                                                       1 & 1 & 0 & 0 & 0 \cr 
                                                       1 & 1 & 1 & 0 & 0 \cr 
                                                       0 & 1 & 1 & 1 & 0 \cr 
                                                       0 & 1 & 1 & 0 & 1 \cr }}\right]$$
 
However, surgering the chords in $D_1$ yields a string link whose closure has one
component, and surgering the chords in $D_2$ yields a string link whose closure has 3
components.  Therefore, $C(D_1) = 1$, but $C(D_2) = 0$.

However, we can at least give a sufficient condition for $C(D)$ to be trivial which
depends only on the adjacency matrix.

\begin{prop} \label{P:suff}
Let D be a string link chord diagram of degree k on n components.  If either of the
following conditions holds, then C(D) = 0.
\begin{itemize}
     \item  $det(\Gamma(D)) = 0$.
     \item  $k+n$ is even.
\end{itemize}
\end{prop}
{\sc Proof:}  The first condition means that, in the normal form, there will be an isolated
chord in the caravan, so $C(D)$ will be trivial.  The second condition comes from observing
that every surgery of a chord in $D$ the number of components in the closure
of $D$ by 1 (either adding or removing a component).  Since we begin with $n$
components, the number of components after all the chords are surgered is congruent to
$k+n$ modulo 2.  So $C(D)$ can be non-trivial only if $k+n \equiv 1$ mod 2.  Hence, if $k+n$
is even, then $C(D) = 0$. $\Box$

Now we will consider the (framed) Homfly weight system for string links.  The Homfly
invariant is the Laurent polynomial $P(l,m) \in {\bf Z}[l^{\pm 1},m^{\pm 1}]$ defined by
the following skein relations \cite{li} ($L^+$ is the result of adding a positive kink to
the link $L$):
$$P(L_+)-P(L_-) = mP(L_0)$$
$$P(L^+) = lP(L)$$
$$P(L \cup O) = \frac{l-l^{-1}}{m}P(L)$$
$$P(O) = 1$$
If we make the substitutions $m = e^{ax/2}-e^{-ax/2}$ and $l =
e^{abx/2}$, and expand the resulting power series, we transform the
Homfly polynomial into a power series in $x$, whose coefficients
are finite type invariants (of regular isotopy).  These invariants
give rise to regular weight systems which we can collect together
as the Homfly regular weight system $H$.  The skein relations
above give rise to the following relations for $H$, by looking at
the first terms of the power series (as before, $D_v$ is the
result of surgering the chord $v$ in $D$):
$$H(D) = aH(D_v)$$
$$H(D \cup O) = bH(D)$$
$$H(O) = 1$$
So if $D$ is an unlink of $k$ components, $H(D) = b^{k-1}$.  Since
the first of these relations is almost the same as for the Conway
weight system $C$, the same argument shows that $H$ satisfies the
2-term relations.  We can now consider string link diagrams with two components (as with
the Conway weight system, the Homfly weight system is not necessarily determined by the
adjacency matrix of the intersection graph for chord diagrams on more than two components).

\begin{thm} \label{T:homfly2}
For any chord diagram D of degree k on a string link with 2 components, let
$r = rank(\Gamma(D))$.  If r is odd, then $H(D) = a^kb^{k-r}$; if r is even, then $H(D) =
a^kb^{k-r+1}$.
\end{thm}
{\sc Proof:}  As with Theorem \ref{T:conway2}, it suffices to show
that the weight systems agree on diagrams in normal form.  Let $D$ be the diagram
with $m_1$ one-humped camels and $m_2$ two-humped camels on component 1, and $l$ chords
between the two components ($l=0,1,2$) (so the degree of $D$ is
$k=l+m_1+2m_2$).  As before, $adj(\Gamma(D)) \cong [1]^l \oplus
[0]^{m_1} \oplus \left[{\matrix{0 & 1 \cr 1 & 0}}\right]^{m_2}$, so
the rank is $r=l+2m_2$.  $r$ is even when $l=0,2$ and odd when $l=1$.

On the other hand, if we surger all the chords (each time multiplying $H$ by
$a$), the resulting link has $m_1 + 2$ components when $l=0,2$ (and $r$ is
even), and $m_1 + 1$ components when $l=1$ (and $r$ is odd).  So when $r$ is
even, $H(D) = a^kb^{m_1 + 1} = a^kb^{k-r+1}$, and when $r$ is odd, $H(D) =
a^kb^{m_1} = a^kb^{k-r}$. $\Box$\\
\\
\noindent {\sc Remark:}  If we let $a = 1$ and $b = 0$, then $H(D) = C(D)$.  According to
Theorem \ref{T:homfly2}, if $D$ is a string link chord diagram with two components, we will
have $H(D) = 0$ in this case {\it unless} $r = k$ is odd (in which case $H(D) = 0^0 = 1$). 
But this is exactly the case when $det(\Gamma(D)) = 1$ (since the matrix has full rank) and
$rank(\Gamma(D)) \equiv 1$ mod 2, so $C(D) = 1$.  So the formulas for the Conway and Homfly
weight systems agree in this case.

\section{Further Questions and Problems} \label{S:questions}

There are several obvious questions and directions for further research.

\begin{quest}
To what extent are the Conway and Homfly weight systems for string links with more than two
components determined by the intersection graph?
\end{quest}

Figure \ref{F:counterexample} gives an example of chord diagrams on 4 components where the
adjacency matrix of the intersection graph is insufficient to determine the Conway and
Homfly weight systems.  However, this example seems to depend strongly on there being 4
components - it is not clear how to find connected examples with, say, 5 components where
the adjacency matrix is insufficient.  And certainly the adjacency matrix (as shown by the
example in Figure \ref{F:counterexample} does not contain all the information of the
intersection graph - in particular, the information contained in the vertex labels.  This
leads to our next question:

\begin{quest}
Can we find invariants of intersection graphs other than the adjacency matrix which give
rise to weight systems?
\end{quest}

\begin{quest}
Can we find other invariants of string links whose weight systems can be computed via
intersection graphs?
\end{quest}

The author has looked at Milnor's homotopy invariants for string links \cite{me5}, but
there are many others.

\begin{quest}
To what extent does the intersection graph determine the chord diagram?
\end{quest}

This question is still open for knots, as well, though some progress has been made 
\cite{cdl, me}.  Some work has been done by the author \cite{me3} for the case when the
intersection graph is a tree, but this is only a bare beginning.

\begin{quest}
Can we generalize these constructions from string links to links?
\end{quest}

We should be able to define intersection graphs for links as a quotient space of the
intersection graphs for string links.  It is yet to be seen whether this is useful.

\small

\normalsize

\end{document}